\documentclass[12pt]{article}

\usepackage{a4,graphicx,amsmath,amsfonts,amssymb,mathrsfs,bbm,color}
\usepackage{caption,subcaption}

\newcommand{\R}{\mathbbm{R}}

\newcommand{\N}{\mathbbm{N}}
\newcommand{\htwo}{\mathcal{H}_2}
\newcommand{\ltworandom}{\mathcal{L}^2(\mathcal{M},\rho)}

\begin{document}


\begin{center}
{\bf \Large Energy-based model order reduction \\[0.1ex] for linear stochastic Galerkin systems \\[0.7ex] of second order}

\vspace{10mm}

{\large Roland~Pulch} \\[1ex] 
{\small Institute of Mathematics and Computer Science, 
Universit\"at Greifswald, \\ 
Walther-Rathenau-Stra{\ss}e~47, 17489 Greifswald, Germany. \\
Email: {\tt roland.pulch@uni-greifswald.de}}
\end{center}

\bigskip


\begin{center}
{\bf Abstract}

\begin{tabular}{p{12.5cm}}
  We consider a second-order linear system of ordinary differential 
  equations (ODEs) including random variables. 
  A stochastic Galerkin method yields a larger deterministic linear 
  system of ODEs. 
  We apply a model order reduction (MOR) of this high-dimensional 
  linear dynamical system, where its internal energy represents a  
  quadratic quantity of interest. 
  We investigate the properties of this MOR with respect to 
  stability, passivity, and energy dissipation. 
  Numerical results are shown for a system modelling a 
  mass-spring-damper configuration. 
\end{tabular}
\end{center}


\section{Introduction}
Mathematical models typically include physical parameters or other 
parameters, which are often affected by uncertainties. 
A well-known approach is to change the parameters into random variables 
to address their variability, see~\cite{sullivan:book}. 
Consequently, an uncertainty quantification (UQ) can be performed. 

We study second-order linear systems of ordinary differential equations 
(ODEs), which contain independent random variables.
Each second-order linear system of ODEs together with its  
internal energy is equivalent to a first-order port-Hamil\-to\-nian (pH) 
system, where the Hamiltonian function represents the internal energy, 
see~\cite{beattie-etal}. 
We use a stochastic Galerkin technique, see~\cite{sullivan:book}, 
which produces a larger deterministic system of second-order linear ODEs. 
The stochastic Galerkin projection is structure-preserving. 
Hence the stochastic Galerkin system also features an internal energy, 
which represents a quadratic output of the linear dynamical system. 

Since the stochastic Galerkin system is high-dimensional, 
we employ a model order reduction (MOR), see~\cite{antoulas}, 
to diminish the dimensionality. 
MOR of linear stochastic Galerkin systems with linear outputs was applied 
in~\cite{freitas-etal,pulch-matcom,pulch-jmi}, for example. 
Now we investigate an MOR, where the internal energy is defined as 
the quantity of interest (QoI).
In~\cite{benner-goyal-duff}, a balanced truncation technique was derived 
to reduce a first-order linear system of ODEs with quadratic output. 
We apply the balanced truncation to the canonical first-order system, 
which is equivalent to the second-order stochastic Galerkin system. 

A reduced system of ODEs exhibits a quadratic output, 
which approximates the underlying internal energy. 
An a posteriori error bound is computable for the quadratic output 
in any MOR method, provided that the systems are asymptotically stable. 
Moreover, we study the properties of the reduced systems 
with respect to dissipation inequalities and passivity. 
A concept to measure a loss of passivity is introduced. 
Finally, we present results of numerical experiments 
using a model of a mass-spring-damper system. 


\section{Problem Definition}
A stochastic modelling is applied to second-order linear dynamical systems, 
which include uncertain parameters. 

\subsection{Second-order linear ODEs including Parameters}
We consider second-order linear systems of ODEs in the form
\begin{equation} \label{eq:odes} 
M(\mu) \ddot{p} + D(\mu) \dot{p} + K(\mu) p = B(\mu) u , 
\end{equation}
where the symmetric matrices $M,D,K \in \R^{n \times n}$ and 
the matrix $B \in \R^{n \times n_{\rm in}}$ depend on parameters 
$\mu \in \mathcal{M} \subseteq \R^q$. 
Input signals $u : [0,\infty) \rightarrow \R^{n_{\rm in}}$ 
are supplied to the system. 
The state variables $p : [0,\infty) \times \mathcal{M} \rightarrow \R^n$ 
depend on time as well as the parameters. 
We assume that the matrices~$M$ and~$K$ are positive definite and 
the matrix~$D$ is positive definite or semi-definite 
for all $\mu \in \mathcal{M}$. 
It follows that each linear dynamical system~(\ref{eq:odes}) is 
Lyapunov stable. 
A positive definite matrix~$D$ is sufficient for the asymptotic stability 
of a system~(\ref{eq:odes}), see~\cite{inman}.

The linear dynamical system~(\ref{eq:odes}) features an internal energy 
\begin{equation} \label{eq:energy}
V(p,\dot{p},\mu) = \tfrac{1}{2} 
\left( \dot{p}^\top M(\mu) \dot{p} + p^\top K(\mu) p \right) , 
\end{equation}
which represents the sum of kinetic energy and potential energy. 
In~\cite{beattie-etal}, it is shown that a second-order linear system 
of ODEs, which satisfies the above assumptions on the definiteness of the 
matrices, is equivalent to a first-order pH system. 
Consequently, the internal energy~(\ref{eq:energy}) is identical to 
the Hamiltonian function of the pH system. 

\subsection{Stochastic Modelling and Polynomial Chaos \\ Expansions}
Often the parameters are affected by uncertainties. 
In UQ, a typical approach consists in replacing the parameters 
by random variables, see~\cite{sullivan:book}. 
Thus we substitute the parameters in the system~(\ref{eq:odes}) by 
independent random variables $\mu : \Omega \rightarrow \mathcal{M}$, 
$\omega \mapsto (\mu_1(\omega),\ldots,\mu_q(\omega))$
on a probability space $(\Omega,\mathcal{A},\mathcal{P})$. 
We use traditional probability distributions for each parameter 
like uniform distribution, beta distribution, Gaussian distribution, etc. 
Let a joint probability density function $\rho : \mathcal{M} \rightarrow \R$ 
be given.  
A measurable function $f : \mathcal{M} \rightarrow \R$ exhibits 
the expected value
\begin{equation} \label{eq:expected-value}
\mathbb{E} [f] = \int_{\Omega} f(\mu(\omega)) \; 
\mathrm{d}\mathcal{P}(\omega) 
= \int_{\mathcal{M}} f(\mu) \, \rho(\mu) \; \mathrm{d}\mu . 
\end{equation}
The expected value~(\ref{eq:expected-value}) implies an inner product 
$\langle f,g \rangle = \mathbb{E}[fg]$ for two square-integrable 
functions~$f,g$. 
We denote the associated Hilbert space by $\ltworandom$. 

Let an orthonormal basis $(\Phi_i)_{i \in \N}$ be given, 
which consists of polynomials $\Phi_i : \mathcal{M} \rightarrow \R$. 
It holds that $\langle \Phi_i , \Phi_j \rangle = \delta_{ij}$ 
with the Kronecker-delta. 
The number of basis polynomials up to a total degree~$d$ is 
$s = \frac{(d+q)!}{d!q!}$. 
This number becomes high for larger~$q$ even if $d$ is moderate, 
say $d \le 5$. 

This orthonormal basis allows for expansions in the so-called 
polynomial chaos (PC), see~\cite{sullivan:book}. 
A function $f \in \ltworandom$ can be represented as 
a PC expansion
\begin{equation} \label{eq:pce}
f(\mu) = \sum_{i=1}^{\infty} f_i \Phi_i(\mu) 
\qquad \mbox{with} \qquad
f_i = \langle f, \Phi_i \rangle . 
\end{equation} 
We apply this expansion to the state variables in~(\ref{eq:odes}) 
separately for each component $p_1,\ldots,p_n$ and each 
time point $t \ge 0$.  

\subsection{Stochastic Galerkin System}
Using the expansion~(\ref{eq:pce}) for the state variables, 
we arrange a finite sum with $s$~terms including a priori unknown 
approximations of the coefficients.  
Inserting the finite sum into~(\ref{eq:odes}) generates a residual. 
The Galerkin approach requires that this residual is orthogonal to 
the subspace $\{\Phi_1,\ldots,\Phi_s\}$ spanned by the basis polynomials.  
The orthogonality is defined using the inner product of the 
Hilbert space $\ltworandom$. 
The stochastic Galerkin projection yields a deterministic 
second-order linear system of ODEs
\begin{equation} \label{eq:galerkin} 
\hat{M} \ddot{\hat{p}} + \hat{D} \dot{\hat{p}} + \hat{K} \hat{p} = 
\hat{B} u 
\end{equation}
with larger matrices $\hat{M},\hat{D},\hat{K} \in \R^{ns \times ns}$, 
and $\hat{B} \in \R^{ns \times n_{\rm in}}$. 
The solution of the system is $\hat{p} : [0,\infty) \rightarrow \R^{ns}$ 
with $\hat{p} = (\hat{p}_1^\top,\ldots,\hat{p}_s^\top)^\top$, 
where $\hat{p}_i$ represents an approximation of the exact 
PC coefficients with respect to the $i$th basis polynomial. 
More details on the stochastic Galerkin projection for linear ODEs 
can be found in~\cite{pulch-matcom,pulch2023}, for example. 

The stochastic Galerkin projection is structure-preserving. 
Thus the matrices $\hat{M},\hat{D},\hat{K}$ are symmetric again and also 
inherit the definiteness of the original matrices $M,D,K$. 
The stochastic Galerkin system~(\ref{eq:galerkin}) exhibits the 
internal energy
\begin{equation} \label{eq:energy-galerkin}
\hat{V}(\hat{p},\dot{\hat{p}}) = \tfrac{1}{2} 
\big( \dot{\hat{p}}^\top \hat{M} \dot{\hat{p}} + 
\hat{p}^\top \hat{K} \hat{p} \big) .  
\end{equation}
The linear dynamical system~(\ref{eq:galerkin}) without input 
($u \equiv 0$) satisfies the dissipation property  
\begin{equation} \label{eq:galerkin-dissipation} 
\tfrac{\mathrm{d}}{\mathrm{d}t} \hat{V}(\hat{p},\dot{\hat{p}}) 
= - \dot{\hat{p}}^\top \hat{D} \dot{\hat{p}} \le 0 ,
\end{equation} 
since we assume that the matrix~$\hat{D}$ is positive (semi-)definite. 

Furthermore, the second-order linear system~(\ref{eq:galerkin}) 
has an equivalent linear explicit first-order system
\begin{equation} \label{eq:galerkin-firstorder} 
\begin{pmatrix}
\dot{\hat{v}}_1 \\ \dot{\hat{v}}_2 \\
\end{pmatrix} = 
\begin{pmatrix}
0 & I_n \\
-\hat{M}^{-1} \hat{K} & -\hat{M}^{-1} \hat{D} \\
\end{pmatrix}
\begin{pmatrix}
\hat{v}_1 \\ \hat{v}_2 \\
\end{pmatrix} +
\begin{pmatrix}
0 \\ \hat{M}^{-1} \hat{B} \\
\end{pmatrix} 
u
\end{equation}
with $\hat{v}_1 = \hat{p}$, $\hat{v}_2 = \dot{\hat{p}}$, 
and identity matrix $I_n \in \R^{n \times n}$. 
The internal energy~(\ref{eq:energy-galerkin}) represents a 
quadratic output of~(\ref{eq:galerkin-firstorder}) due to
\begin{equation} \label{eq:energy-galerkin-firstorder}
\hat{V} ( \hat{v}_1 , \hat{v}_2 ) 
= \tfrac{1}{2} 
\begin{pmatrix} \hat{v}_1 \\ \hat{v}_2 \\ \end{pmatrix}^\top
\begin{pmatrix} \hat{K} & 0 \\ 0 & \hat{M} \\ \end{pmatrix} 
\begin{pmatrix} \hat{v}_1 \\ \hat{v}_2 \\ \end{pmatrix} . 
\end{equation}
This relation is shortly written as 
$\hat{V}(\hat{v}) = \frac12 \hat{v}^\top \hat{N} \hat{v}$.


\section{Dissipation Inequality and Passivity}
\label{sec:dissipativity}
Let a linear dynamical system be given in the form $\dot{x} = A x + B u$ 
with $A \in \R^{n \times n}$ and $B \in \R^{n \times n_{\rm in}}$. 
The quadratic output $V = \frac12 x^\top N x$ with $N \in \R^{n \times n}$ 
satisfies the dissipation inequality 
\begin{equation} \label{eq:dissipativity} 
\tfrac{\mathrm{d}}{\mathrm{d}t} x^\top N x \le 
u^\top R u + 2 u^\top S x + x^\top L x 
\end{equation}
with two symmetric matrices~$L \in \R^{n \times n}$, 
$R \in \R^{n_{\rm in} \times n_{\rm in}}$, 
and matrix $S \in \R^{n_{\rm in} \times n}$,   
if and only if the symmetric matrix
\begin{equation} \label{eq:matrix-property}
\left( \begin{array}{c|c}
A^\top N + N A - L & N B - S^\top \\ \hline
B^\top N - S & - R \\
\end{array} \right) 
\end{equation}
is negative definite or semi-definite, see~\cite{willems}.  
We select~$R=0$ and $S=B^\top N$. 
Advantageous is a bound~(\ref{eq:dissipativity}) with $L=0$, 
because this case implies a dissipation inequality
\begin{equation} \label{eq:dissipation-passive} 
\tfrac{\mathrm{d}}{\mathrm{d}t} \tfrac12 x^\top N x \le u^\top y
\end{equation}
including the linear output $y = B^\top N x$, 
as in pH systems. 
Consequently, the linear dynamical system is passive, 
see~\cite{schaft-jeltsema}. 
Usually, the term $u^\top y$ is interpreted as supplied power 
and the term $\tfrac12 x^\top N x$ as internal energy or stored energy. 
Thus we insert $R=0$, $L=0$, $S = B^\top N$ in~(\ref{eq:matrix-property}). 
It follows that the passivity condition~(\ref{eq:dissipation-passive}) 
is satisfied, if and only if the matrix $A^\top N + N A$ 
is negative definite or semi-definite.  


\section{Model Order Reduction}
We perform an MOR of the stochastic Galerkin system, 
where the internal energy represents the QoI. 

\subsection{Model Order Reduction for Linear Systems \\ with Quadratic Output} 
\label{sec:bt}
The full-order model (FOM) is a general first-order linear system of ODEs 
with quadratic output 
\begin{align} \label{eq:fom}
\begin{split}
\dot{x} & = A x + B u \\[-0.5ex]
y & = x^\top N x 
\end{split}
\end{align}
including a symmetric matrix~$N$. 
Let $n$ be the dimension of this system again. 
In~\cite{benner-goyal-duff}, a balanced truncation method was introduced 
for systems of the form~(\ref{eq:fom}). 
This technique requires that the system is asymptotically stable. 
We outline this method. 
The two Lyapunov equations 
\begin{align}
A P + P A^\top + B B^\top & = 0 \label{eq:controllability} \\
A^\top Q + A Q + N P N & = 0    \label{eq:observability}
\end{align} 
are solved successively, which yields the controllability Gramian~$P$ and 
the observability Gramian~$Q$. 
Now symmetric decompositions $P = Z_P Z_P^\top$ and $Q = Z_Q Z_Q^\top$ 
are applied. 
The singular value decomposition (SVD) 
\begin{equation} \label{eq:svd} 
Z_P^\top Z_Q = U \Sigma V^\top 
\end{equation}
yields orthogonal matrices $U,V$ and a diagonal matrix~$\Sigma$, 
which includes the singular values in descending order. 
We choose a reduced dimension $r$. 
Let $U=(U_1,U_2)$, $V=(V_1,V_2)$, 
and $\Sigma = {\rm diag}(\Sigma_1,\Sigma_2)$ 
with $U_1,V_1 \in \R^{n \times r}$, and $\Sigma_1 \in \R^{r \times r}$. 
We obtain projection matrices
$$ V = Z_P U_1 \Sigma_1^{-\frac12} 
\qquad \mbox{and} \qquad 
W = Z_Q V_1 \Sigma_1^{-\frac12} . $$
The reduced-order model (ROM) of dimension~$r$ becomes 
\begin{align} \label{eq:rom}
\begin{split}
\dot{\bar{x}} & = \bar{A} \bar{x} + \bar{B} u \\[-0.5ex]
\bar{y} & = \bar{x}^\top \bar{N} \bar{x}  
\end{split}
\end{align}
with the smaller matrices 
$\bar{A} = W^\top A V$, $\bar{B} = W^\top B$, $\bar{N} = V^\top N V$. 
The linear dynamical system~(\ref{eq:rom}) inherits the asymptotic 
stability of the linear dynamical system~(\ref{eq:fom}) 
in the balanced truncation technique.

Furthermore, an a posteriori error bound can be computed for 
the quadratic output in any MOR method, see~\cite{benner-goyal-duff}. 
We denote the linear dynamical systems~(\ref{eq:fom}) and (\ref{eq:rom}) 
by $H$ and $\bar{H}$, respectively. 
The error of the MOR for the quadratic output is measured in 
the $\htwo$-norm.  
The norm of the system~(\ref{eq:fom}) reads as 
$$ \| H \|_{\htwo} = \sqrt{{\rm trace}(B^\top Q B)} $$
with the observability Gramian satisfying~(\ref{eq:observability}). 
Likewise, we obtain the $\htwo$-norm of the system~(\ref{eq:rom}).
It holds that 
$$ \| y - \bar{y} \|_{\mathcal{L}^{\infty}} \le 
\| H - \bar{H} \|_{\htwo} \| u \otimes u \|_{\mathcal{L}^2} $$ 
using the norms of Lebesgue spaces in time domain.  
The error bound can be computed directly by 
\begin{equation} \label{eq:htwo-error} 
\left\| H - \bar{H} \right\|_{\htwo} = 
\sqrt{ {\rm trace} \big( B^\top Q B + \bar{B}^\top \bar{Q} \bar{B} 
- 2 B^\top Z \bar{B} \big) } . 
\end{equation}
Therein, the matrix~$\bar{Q} \in \R^{r \times r}$ satisfies the Lyapunov 
equation~(\ref{eq:observability}) associated to the ROM~(\ref{eq:rom}). 
The matrix~$Z \in \R^{n \times r}$ solves the Sylvester equation 
\begin{equation} \label{eq:sylvester1} 
A^\top Z + Z \bar{A} + N X \bar{N} = 0 , 
\end{equation} 
while $X \in \R^{n \times r}$ represents the solution of 
the Sylvester equation
\begin{equation} \label{eq:sylvester2} 
A X + X \bar{A}^\top + B \bar{B}^\top = 0 . 
\end{equation}
Lyapunov equations and Sylvester equations can be solved numerically 
either by direct methods or iterative methods. 

\subsection{Application to Stochastic Galerkin System}
The second-order stochastic Galerkin system~(\ref{eq:galerkin}) 
and its internal energy~(\ref{eq:energy-galerkin}) is equivalent to  
the first-order system~(\ref{eq:galerkin-firstorder}) 
with quadratic output~(\ref{eq:energy-galerkin-firstorder}). 
The dissipation analysis of Section~\ref{sec:dissipativity} can be 
applied to 
(\ref{eq:galerkin-firstorder}), (\ref{eq:energy-galerkin-firstorder}). 
We obtain
\begin{equation} \label{eq:galerkin-red-diss} 
\hat{A}^\top \hat{N} + \hat{N} \hat{A} = 
\begin{pmatrix} 
0 & 0 \\ 0 & - 2 \hat{D} \\
\end{pmatrix} . 
\end{equation}
The positive (semi-)definiteness of the matrix~$\hat{D}$ is equivalent 
to the negative definiteness of the 
matrix~(\ref{eq:galerkin-red-diss}). 
Thus the stochastic Galerkin system features the desired 
dissipation inequality~(\ref{eq:dissipation-passive}) 
and thus it is passive. 
This property of the matrix~(\ref{eq:galerkin-red-diss}) 
is related to the counterpart~(\ref{eq:galerkin-dissipation}).

We employ the MOR method from Section~\ref{sec:bt} to the 
high-dimensional system~(\ref{eq:galerkin-firstorder}) 
with quadratic output~(\ref{eq:energy-galerkin-firstorder}). 
The balanced truncation technique preserves the asymptotic stability 
of the FOM, i.e., 
each ROM is asymptotically stable again. 
However, the balanced truncation technique does not preserve the 
passivity with respect to the internal energy, 
as demonstrated by a test example in Section~\ref{sec:numerical-results}. 
Hence the matrix 
\begin{equation} \label{eq:matrix-T} 
\bar{T} := \bar{A}^\top \bar{N} + \bar{N} \bar{A}
\end{equation} 
is not negative (semi-)definite in general. 
Let $\lambda_{\max} > 0$ be the largest eigenvalue of $\bar{T}$. 
A shift of the spectrum via $\bar{T} - \lambda_{\max} I_r$ with 
identity matrix $I_r \in \R^{r \times r}$ yields a 
negative semi-definite matrix. 
Choosing $\bar{R}=0$, $\bar{S} = \bar{B}^\top \bar{N}$, 
$\bar{L} = \lambda_{\max} I_r$ implies the dissipation inequality, 
cf.~(\ref{eq:dissipativity}), 
\begin{equation} \label{eq:dissipation-lambda} 
\tfrac{\mathrm{d}}{\mathrm{d}t} \bar{x}^\top \bar{N} \bar{x} 
\le 2 u^\top \bar{B}^\top \bar{N} \bar{x} + \lambda_{\max} \, x^\top x
= 2 u^\top \bar{B}^\top \bar{N} \bar{x} + \lambda_{\max} \, \| x \|_2^2 . 
\end{equation}
The desired property would be the case of $\lambda_{\max} \le 0$. 
Hence the magnitude of~$\lambda_{\max} > 0$ measures the loss of 
passivity.


\section{Numerical Results}
\label{sec:numerical-results}
\begin{figure}
\begin{center}
\includegraphics[width=40mm]{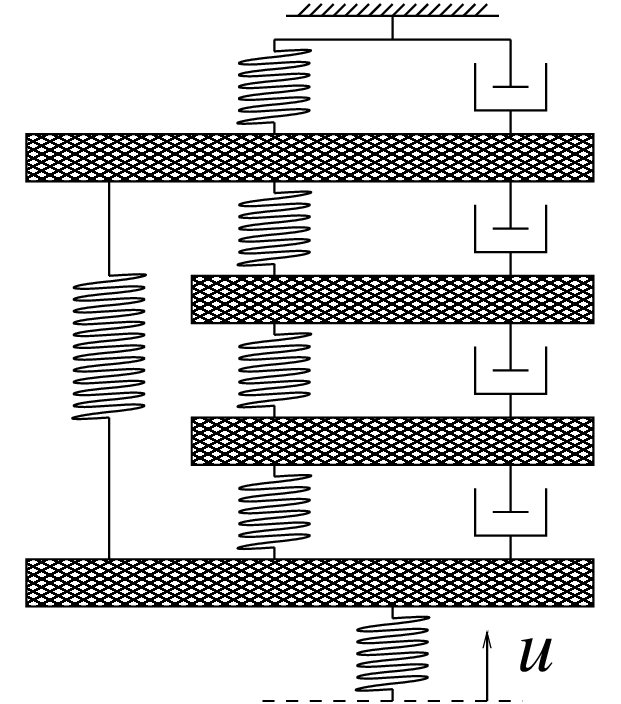}
\end{center}
\caption{Mass-spring-damper configuration.}
\label{fig:mass-spring-damper}
\end{figure}

As test example, we employ a mass-spring-damper system 
from~\cite{lohmann-eid}. 
Figure~\ref{fig:mass-spring-damper} shows the configuration. 
The system contains 4~masses, 6~springs, and 4~dampers, in total
$q=14$ physical parameters.  
A single input~$u$ is supplied by an excitation at the lowest spring. 
This test example was also used in~\cite{pulch-matcom,pulch2023}. 
The mathematical model consists of $n=4$ second-order ODEs~(\ref{eq:odes}). 
The matrices $M,K,D$ are symmetric as well as positive definite for 
all positive parameters. 

In the stochastic modelling, we replace the parameters by random variables 
with independent uniform probability distributions, which vary 
10\% around their mean values. 
Consequently, the PC expansions~(\ref{eq:pce}) include the 
(multivariate) Legendre polynomials. 
We study two cases of total degree: two and three. 
Table~\ref{tab:galerkin-system} demonstrates the properties of the 
resulting second-order stochastic Galerkin systems~(\ref{eq:galerkin}). 
In particular, the sparsity of the system matrices is specified 
by the percentage of non-zero entries. 
The stochastic Galerkin systems are asymptotically stable, 
since the Galerkin projection preserves the definiteness of matrices.   

\begin{table}[b]
\caption{Properties of second-order stochastic Galerkin system.}
\label{tab:galerkin-system}
\begin{center}
\begin{tabular}{@{}cccccc@{}}
\hline
& no. of basis & & non-zero & non-zero & non-zero \\[-0.5ex]
degree & polynomials & dimension & 
entries in~$\hat{M}$ & entries in~$\hat{D}$ & entries in~$\hat{K}$ \\
\hline
2 & 120 & 480 & 0.26\% & 0.69\% & 0.86\% \\
3 & 680 & 2720 & 0.05\% & 0.13\% & 0.17\% \\ 
\hline
\end{tabular}
\end{center}
\end{table}

\begin{figure}
\centering
\begin{subfigure}[b]{0.45\textwidth}
         \centering
         \includegraphics[width=\textwidth]{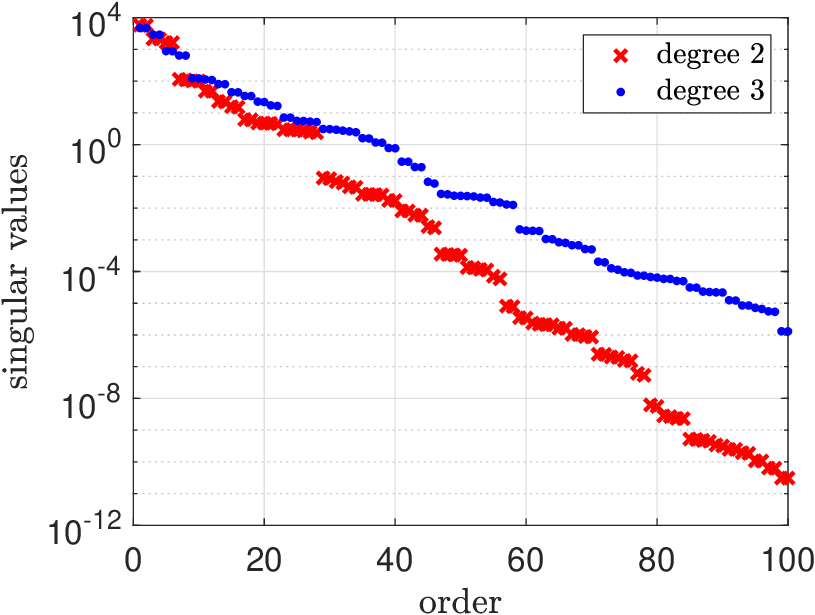}
         \caption{singular values}
     \end{subfigure}
     \hspace{5mm}
     \begin{subfigure}[b]{0.45\textwidth}
         \centering
         \includegraphics[width=\textwidth]{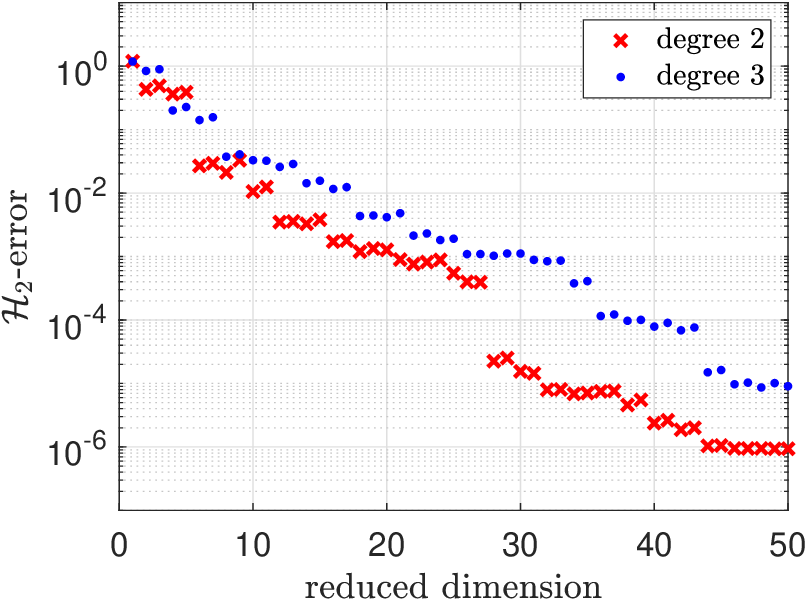}
         \caption{relative $\htwo$-error of MOR}
     \end{subfigure}
\caption{Singular values of Hankel-type and relative $\htwo$-error of 
internal energy in balanced truncation technique.}
\label{fig:mor}
\end{figure}

Now we perform an MOR of the equivalent first-order 
system~(\ref{eq:galerkin-firstorder}) 
with quadratic output~(\ref{eq:energy-galerkin-firstorder}) 
using the balanced truncation technique from Section~\ref{sec:bt}. 
We solve the Lyapunov equations 
(\ref{eq:controllability}), (\ref{eq:observability}) 
and the Sylvester equations~(\ref{eq:sylvester1}), (\ref{eq:sylvester2}) 
by direct methods of numerical linear algebra. 
Figure~\ref{fig:mor} (a) depicts the Hankel-type 
singular values of the SVD~(\ref{eq:svd}), 
which rapidly decay to zero. 
We compute the ROMs~(\ref{eq:rom}) of dimension $r=1,\ldots,100$. 
The error of the MOR is measured in the relative $\htwo$-norm, 
i.e., $\| H - \bar{H} \|_{\htwo} / \| H \|_{\htwo}$, 
see~(\ref{eq:htwo-error}). 
The relative errors are shown for $r \le 50$ 
in Figure~\ref{fig:mor} (b). 
We observe that a high accuracy is achieved already for relatively  
small reduced dimensions. 

\begin{figure} 
\centering
\begin{subfigure}[b]{0.45\textwidth}
         \centering
         \includegraphics[width=\textwidth]{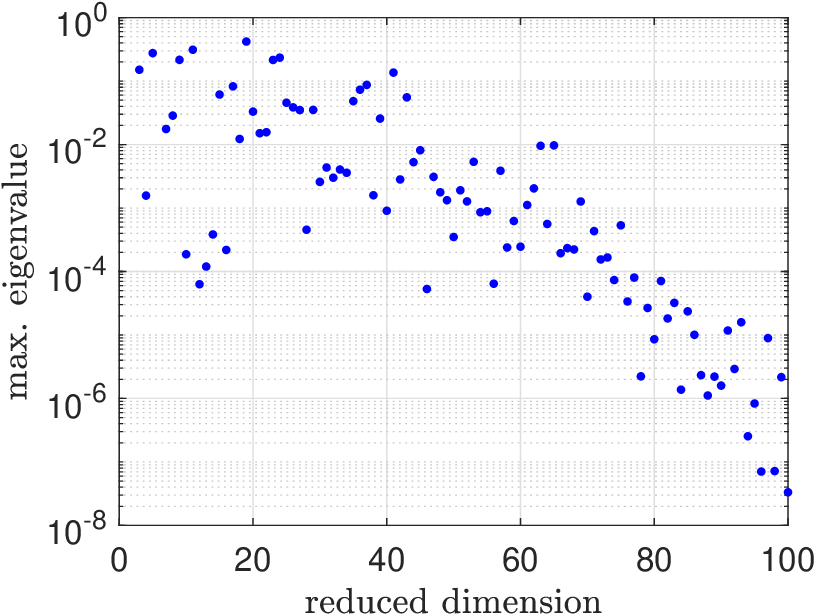}
         \caption{degree two}
     \end{subfigure}
     \hspace{5mm}
     \begin{subfigure}[b]{0.45\textwidth}
         \centering
         \includegraphics[width=\textwidth]{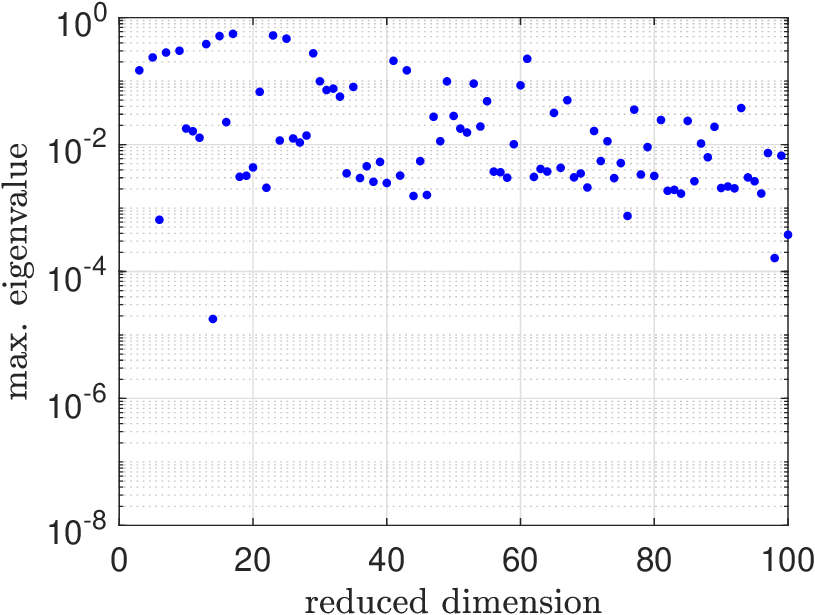}
         \caption{degree three}
     \end{subfigure}
\caption{Maximum eigenvalue of matrices~(\ref{eq:matrix-T}) using 
balanced truncation in the case of polynomial degree two and three.}
\label{fig:eigenvalues}
\end{figure}

Furthermore, we examine the dissipation properties of the 
ROMs~(\ref{eq:rom}), as described in Section~\ref{sec:dissipativity}.
All reduced systems loose the passivity, 
because their matrices~(\ref{eq:matrix-T}) 
are not negative (semi-)definite. 
The maximum eigenvalues of the matrices are illustrated by 
Figure~\ref{fig:eigenvalues}. 
The maxima tend to zero for increasing reduced dimension. 
Yet the decay becomes slower for larger total polynomial degree. 
It follows that the dissipation inequality~(\ref{eq:dissipation-lambda}) 
is valid for a small eigenvalue~$\lambda_{\max}$. 

\begin{figure}
\centering
\begin{subfigure}[b]{0.45\textwidth}
         \centering
         \includegraphics[width=\textwidth]{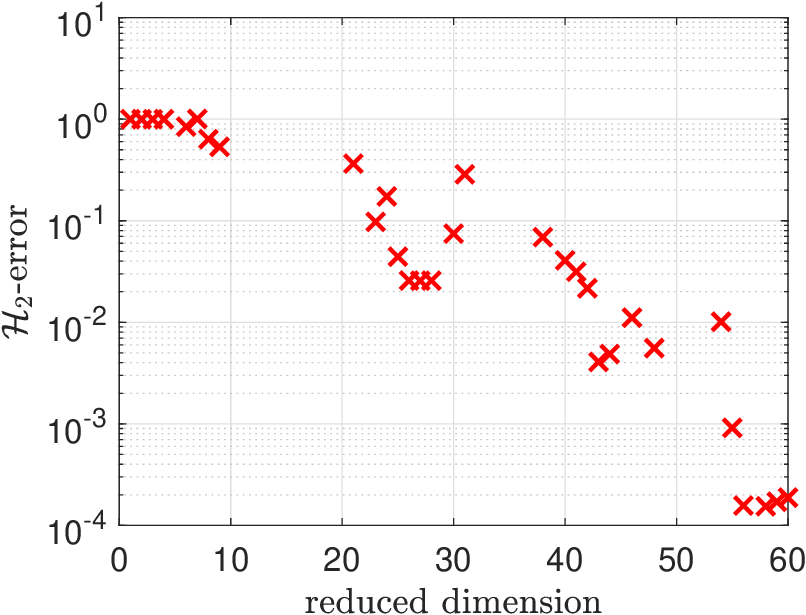}
         \caption{relative $\htwo$-error of MOR}
     \end{subfigure}
     \hspace{5mm}
     \begin{subfigure}[b]{0.45\textwidth}
         \centering
         \includegraphics[width=\textwidth]{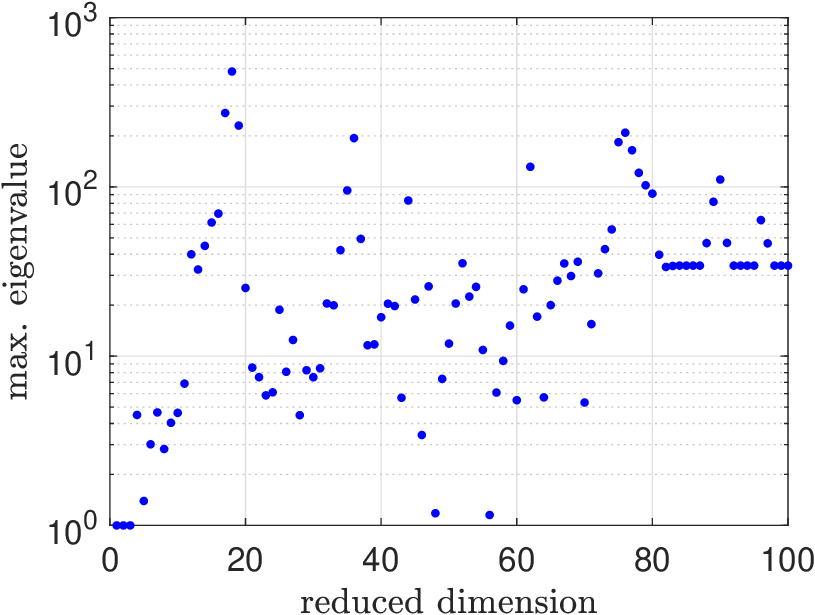}
         \caption{maximum eigenvalue}
     \end{subfigure}
\caption{Relative $\htwo$-error of internal energy in MOR 
and maximum eigenvalue of matrices~(\ref{eq:matrix-T})
in Arnoldi method for polynomial degree two.}
\label{fig:arnoldi}
\end{figure}

Finally, we present a comparison. 
We reduce the stochastic Galerkin system~(\ref{eq:galerkin-firstorder}) 
for polynomial degree two by the Arnoldi method,  
which is a specific Krylov subspace technique, see~\cite{antoulas}.  
This scheme is a Galerkin-type MOR method, i.e., the 
projection matrices satisfy $V=W$. 
However, the asymptotic stability may be lost in this technique. 
The Arnoldi method does not include an information about 
a definition of a QoI.  
We use a single (real) expansion point $\omega=1$ in the complex 
frequency domain, 
because other real-valued choices $\omega = 10^k$ with 
$k \in \{-2,1,1,2\}$ cause worse approximations. 
Figure~\ref{fig:arnoldi} (a) depicts the 
relative $\htwo$-error of the internal energy for the ROMs 
of dimension $r \le 60$. 
Higher reduced dimensions produce larger errors due to 
an accumulation of round-off errors in the orthogonalisation, 
which is a well-known effect in the Arnoldi algorithm.   
If an ROM~(\ref{eq:rom}) is unstable, then the error is not computable 
and thus omitted. 
As expected, the accuracy of the Arnoldi method is not as good 
as the accuracy of the balanced truncation. 
Again the passivity is lost in all ROMs. 
Figure~\ref{fig:arnoldi} (a) shows 
the maximum eigenvalue of the matrices~(\ref{eq:matrix-T}).  
We observe that these positive maxima do not decay, 
even though small errors are achieved for reduced dimensions 
$55 \le r \le 60$. 

\section{Conclusions}
We applied a stochastic Galerkin projection to a second-order linear 
system of ODEs including random variables. 
The high-di\-men\-sional stochastic Galerkin system owns an internal 
energy as quadratic output. 
We performed an MOR of an equivalent first-order system of ODEs, 
where the used balanced truncation method is specialised to approximate 
a quadratic output. 
However, the passivity of the dynamical systems with respect to the 
internal energy may be lost in this reduction. 
We proposed a concept to quantify the discrepancy of a non-passive 
dynamical system to the passive case. 
Numerical results of a test example demonstrated that this 
discrepancy measure tends to zero for increasing reduced dimensions 
in the balanced truncation method. 

\clearpage


%
\end{document}